\documentclass{ifacmtg}
\usepackage[dcucite]{harvard}

\usepackage[dvips]{graphicx}
\usepackage{makeidx}
\usepackage{amsmath}
\usepackage{amsfonts}

\newcommand {\etavec}{\boldsymbol{\eta}}

\newcommand {\xivec}{{\boldsymbol{\xi}}}

\newcommand {\thetavec}{{\boldsymbol{\theta}}}

\newfont{\pseudocode}{cmtt10}

\newtheorem{assume}{Assumption}

\newcommand{\bfx}{\mathbf{x}}
\newcommand{\bff}{\mathbf{f}}

\newcommand{\bfg}{\mathbf{g}}
\newcommand{\const}{\mathrm{const}}

\newcommand{\dist}{\mathrm{dist}}

\newcommand{\im}{\mathcal{I}m}

\newcommand{\dpsi}{{\dot\psi}}
\newcommand{\pd}{{\partial}}

\newcommand{\Real}{\mathbb{R}}
\newcommand{\Natural}{\mathbb{N}}

\makeindex

\begin{document}



\runauthor{Ivan Tyukin, Cees van Leeuwen}
\begin{frontmatter}
\title{Adaptation and nonlinear parametrization: nonlinear dynamics prospective}
\author[Wako-shi]{Ivan Tyukin}
\author[Wako-shi]{Cees van Leeuwen}
\address[Wako-shi]{Lab. for Perceptual Dynamics,
RIKEN Brain Science Institute, 2-1, Hirosawa, Wako-shi, Saitama,
Japan}

\begin{abstract}
We consider adaptive control problem in presence of nonlinear
parametrization of uncertainties in the model. It is shown that
despite traditional approaches require for domination in the
control loop  during adaptation, it is not often necessary to use
such energy inefficient compensators it in wide range of
applications. In particular, we show that recently introduced
adaptive control algorithms in finite form which are applicable to
monotonic parameterized systems  can be extended to general smooth
non-monotonic parametrization. These schemes do not require any
damping or domination in control inputs.
\end{abstract}

\begin{keyword}
adaptive systems, nonlinear parameterization, finite form
algorithms, convergence, nonlinear persistent excitation
\end{keyword}
\end{frontmatter}

\section{Introduction}

The standard techniques to address nonlinear parameterization of
the uncertainty in adaptive control literature involve applying
domination functions (see, for instance,
\cite{MarinoTomei93_nonlin,Lin,Lin_2002_smooth}) or damping of the
unknown nonlinearity \cite{Annaswamy99}. Both operate on the
principle of invoking auxiliary crutches ensuring Lyapunov
stability of the extended system. These methods yield effective
analytical solutions to the problem of adaptation in presence of
nonlinear parameterized uncertainties.

Yet, in a wide range of applications in control and reverse
bio-engineering, physics and biology, it is desirable to refrain
from domination or damping. The issue of non-dominating adaptation
is critical in control engineering applications where systematic
overshooting in control results in fast wearing off of the
actuators, waste of energy and undesired chattering. A typical
example is traction/braking control under unknown tire-road
conditions.  Tire-road conditions enter into the equations
describing the slip dynamics as uncertainties that are highly
nonlinear in their parameters \cite{Pacejka91,Canudas_1999}. The
issue is how under these conditions the slip could be effectively
controlled during braking without applying unnecessary large
braking/trackgin torques.

In bio-engineering, the motivation to use non-conventional
adaptation mechanisms is even more stronger. In these systems
physiological considerations motivate local adaptation, at the
level of the single node \cite{Brenner_2000},\cite{Webster_2002}
(visual system, cognitive processing e.t.c.). Even simplest
mathematical models of these nodes usually are nonlinearly
parameterized \cite{Abbott_2001}. Due to the necessity to respond
sensitively to different stimulation, use of external domination
is not desirable. On the other hand, because the number of
interacting nodes is large, the adaptation, should be
non-dominating in order to avoid inefficient consumption of
energy.

From control-theoretic prospective, perhaps, the most challenging
task is to find the principles, which will allow to control
behavior of nonlinear systems by gentle and non-dominating
parametric adaptation. Available non-dominating solutions,
however, are either local \cite{Karsenti_1996}, or assume
monotonic parametrization of uncertainties
\cite{t_fin_formsA&T2003,tpt2003_tac}. Therefore, new approaches
are needed in order to design adaptive algorithms capable of
ensuring asymptotic reaching of the control goal for large class
of nonlinear parameterized systems

Finding gentle solutions to control problems with nonlinear
uncertainties will involve a change in design methodology.
Conventional design methods in adaptive control theory often favor
Lyapunov-based methodology. For many decades Lyapunov methods were
successfully applied in design and analysis of nonlinear and
adaptive systems (see, for example,
\cite{Isidory,Sastry_1999,Narendra89,Fradkov99}). Furthermore, as
a rule of thumb, Lyapunov stability serves as a measure of
acceptable performance in mechanics and engineering simply for its
formulation guaranteeing small deviations from the equilibrium in
case of small perturbations in initial conditions.

When adaptive control is required, however, deviations of the
unknown parameters are likely to be large in order for the
adaptation to make practical sense. Therefore, requirement of
Lyapunov stability for systems with general parametrization does
not seem to be necessary unless it provides certain advantages in
specific applications in addition to mere asymptotic reaching of
the control goal.


On the other hand, adaptation processes in many physical and
biological phenomena are far from being stable in conventional
sense. Examples of such non-stable adaptation include (but are not
limited to) high-frequency oscillation of the eye (tremor) that
allows us to see static pictures better
\cite{Ditchburn_1952,Martinez_2004},  perception of the ambiguous
figures (perceptual switches) \cite{Ito_2003}, evolution in social
and ecological systems \cite{Sole_1999,Bak_1993}. This motivates
us to seek for replacement of conventional Lyapunov
stability-based methodology for dealing with problems of
adaptation in nonlinear systems.

Doing so would enable, in principle, a specter of entirely new
applications and associated problems is emerging. Theoretical
findings in physics and biology such as intermittent
synchronization
\cite{Gauthier_1996,Kaneko_and_Tsuda,Kaneko94,Kaneko_2003},
homeostasis in the living cells stability (see \cite{Moreau_2003}
and references therein) and, last but not least, {\it
self-organized criticality} - phenomena often observed in the
earthquakes \cite{Bak_2002}, in the neuronal activity of the human
brain \cite{Beggs_2003} which can be modeled by the arrays of
phase oscillators with nonlinear driving \cite{Corral_1995} show
the importance of systems on the boundary of stability.
Understanding the principles of learning and adaptation in such
systems could eventually benefit greatly from the principles of
adaptive control.

We start with brief analysis of design strategies for adaptive
systems with nonlinear parametrization. We show that that
certainty-equivalence adaptive control inevitably requires damping
or domination functions to guarantee global stability in presence
of general nonlinear parametrization. In order to avoid these
difficulties we replace requirement of Lyapunov stability of the
system with mere reaching a neighborhood of the goal manifold.
After the new control goal is defined we analyze a class of
nonlinear systems where uncertainties are given by one-dimensional
nonlinear parameterized function. For this class of systems we
provide adaptive control algorithms capable of steering the system
state to a small neighborhood of the target manifold. The control
algorithm should not involve neither domination function, nor
should it require for additional damping unknown nonlinearity.
After this step is accomplished we extend our approach to
multidimensional parameterizations.

Throughout the paper we will use the following notations. Symbol
$\bfx(t,\bfx_0,t_0)$ denotes solution of a system of differential
equations starting at the point $\bfx_0$ at time instant $t_0$;
symbol $C^r$ denotes the space of $r$ times differentiable
functions; symbol $\Real$ stand for the space of reals; $\Real_+$
defines nonnegative real numbers, symbol $\im$ denotes image of
the map. We say that $\nu: \Real_+\rightarrow \Real$  belongs to
$L_2$ iff $L_2(\nu)=\int_0^{\infty}\nu^2(\tau)d\tau<\infty$. The
value $\|\nu\|_2=\sqrt{L_2(\nu)}$ stands for the $L_2$ norm of
$\nu(t)$. Function $\nu: \Real_{+}\rightarrow \Real$ belongs to
and $L_\infty$ iff
$L_{\infty}(\nu)=\sup_{t\geq0}\|\nu(t)\|<\infty$, where
$\|\cdot\|$ is the Euclidean norm. The value of
$\|\nu\|_\infty=L_{\infty}(\nu)$ stands for the $L_{\infty}$ norm
of $\nu(t)$. Symbol $\mathcal{U}_\epsilon(\bfx)$ denotes the set
of all $\bfx': \|\bfx-\bfx'\|\leq\epsilon$. Let
$\mathcal{X}\subset\Real^n$, distance
$\dist(\bfx,\mathcal{X})=\inf_{\bfx'\in
\mathcal{X}}\|\bfx-\bfx'\|$. Symbol
$\mathcal{U}_\epsilon(\mathcal{X})$ denotes the set $\{\bfx'\in
\Real^n| \ \bfx': \ \dist(\bfx',\mathcal{X})\leq \epsilon \}$.
Symbol $L_\bff \psi$ stands for the Lie-derivative of function
$\psi(\bfx)$ w.r.t. the vector field $\bff(\bfx)$.

The paper is organized as follows. Section 2 contain preliminary
notations and statement of the problem. In Section 3 we provide
main results of the paper. Section 4 concludes the paper.

\section{Problem Formulation and Preliminaries}

Let us consider standard adaptive control problem where the
uncertain system be given as follows
\begin{equation}\label{eq:intro:plant}
\dot{\bfx}=\bff(\bfx,\thetavec)+\bfg(\bfx)u, \ \bfx_0\in
\Omega_\bfx
\end{equation}
where $\bfx\in \Real^n$ is state vector,
$\bff:\Real^n\times\Real^d\rightarrow \Real^n$,
$\bfg:\Real^{n}\rightarrow\Real^n$ are $C^1$-smooth vector-fields,
$\thetavec\in \Omega_\theta\subset \Real^d$ is vector of
parameters, $u$ is control input, and $\Omega_\bfx\subset\Real^n$
is the set of initial conditions $\bfx_0$. Functions $\bff$,
$\bfg$ are known, vector of parameters $\thetavec$ is assumed to
be unknown a-priori. We also assume that $\bfg(\bfx)$ is bounded
w.r.t. $\bfx$. Let the control goal be to reach asymptotically a
neighborhood of the following manifold given implicitly by
\begin{equation}\label{eq:intro:goal}
\psi(\bfx)=0, \ \psi\in C^2
\end{equation}
In addition to (\ref{eq:intro:goal}) we will require that
\begin{equation}\label{eq:intro:psi_assume}
\psi(\bfx(t))\in L_\infty \Rightarrow \bfx(t)\in L_\infty
\end{equation}
This requirement ensures that any bounded deviation from the goal
manifold does not result in unbounded growth of the norm
$\|\bfx\|$. In addition we assume that
$|L_\bfg(\bfx)\psi(\bfx))|\geq \delta_{\psi}>0$.

Let us select the class of admissible control functions which,
allowing for the goal to be reached, can compromise between
performance and domination issues. The most natural way would be
to define this class on the ground of {\it certainty-equivalence
principle}.
In particular, consider control function
\begin{equation}\label{eq:intro:control}
\begin{split}
u(\bfx,\hat{\thetavec})&=(L_\bfg\psi(\bfx))^{-1}(-
L_\bff(\bfx,\hat{\thetavec})\psi(\bfx) -\\
& \varphi(\psi)+\upsilon(t))
\end{split}
\end{equation}
which transforms the system (\ref{eq:intro:plant}) into the
following error model
\begin{equation}\label{eq:intro:error}
\dpsi=f(\bfx,\thetavec)-f(\bfx,\hat{\thetavec})-\varphi(\psi)+\upsilon(t),
\end{equation}
where $f(\bfx,\thetavec)=L_\bff(\bfx,\thetavec)\psi(\bfx)$,
$\varphi\in \mathcal{C}_\varphi \subset C^0:\Real\rightarrow\Real$
and $\varphi(\psi)\psi>0 \ \forall \ \psi\neq 0$, $\upsilon, d:
\Real_+\rightarrow \Real$. The desired dynamics (performance) of
the system is given by equation $\dpsi=-\varphi(\psi)$, and
$\upsilon(t)$ stands for auxiliary control or disturbances
depending on the context.

\begin{defn} Adaptive control law (\ref{eq:intro:control})
 is called {\it
non-dominating} in class $\mathcal{C}_\varphi$ if  for any
$\epsilon>0$ there exists such function
$\hat{\thetavec}(\bfx,t,\delta(\epsilon))\in\Omega_{\theta}$,
$\delta(\epsilon):\Real_+\rightarrow\Real$, and time $t^\ast$ that
$|\psi(\bfx(t,\bfx_0,t_0))|<\epsilon$ for any $t\geq t^\ast$,
$\bfx_0\in \Omega_\bfx$, $\thetavec\in\Omega_\theta$, $\varphi\in
\mathcal{C}_\varphi$ and $\upsilon(t)\equiv 0$.
\end{defn}
%
%
In the present study we will restrict class
$\mathcal{C}_{\varphi}$ to the following class of functions:
\begin{equation}\label{eq:intro:class_phi}
\begin{split}
\mathcal{C}_\varphi(k)&=\{\varphi: \Real\rightarrow \Real| \
\varphi\in C^1, \\
&  \varphi(\psi)\psi\geq \psi^{2}k, \ k\in \Real_+\}
\end{split}
\end{equation}
The fact that adaptation is non-dominating in this class of
functions means that for any arbitrary small gain $k>0$ in
feedback $\varphi(\psi)=k\psi$ there exists function
$\hat{\thetavec}(\bfx,t,\delta)$ such that control goal is reached
in finite time. In order to specify desired performance of the
adaptive algorithm itself we require that
$\hat{\thetavec}(\bfx,t,\delta)$ does not change along the
manifold $f(\bfx,\thetavec)-f(\bfx,\hat{\thetavec})=0$ and norm
$\|\thetavec-\hat{\thetavec}\|$ does not grow with time.

In conventional certainty-equivalence adaptive control the problem
of adaptation is usually viewed as the problem of design of
operator $A(\psi,\bfx,\hat{\thetavec},t)$  such that solutions of
\begin{equation}\label{eq:intro:CE_alg}
\dot{\hat{\thetavec}}= A(\psi,\bfx,\hat{\thetavec},t)
\end{equation}
together with (\ref{eq:intro:control}) ensure the goal relation
(\ref{eq:intro:goal}). Function $A(\psi,\bfx,\hat{\thetavec},t)$,
in addition, neither should depend on unknown $\thetavec$, nor it
should require derivative $\dot{\bfx}$.  The standard way to solve
this problem for $\upsilon(t)\equiv 0$ is to design function
$A(\psi,\bfx,\hat{\thetavec},t)$ such that $\psi(\bfx)=0$,
$\hat{\thetavec}=\thetavec$ is stable manifold in Lyapunov sense.

In general nonlinear setup this condition, however, is hardly ever
met for every $\thetavec\in \Omega_\theta$ and $\bfx(t_0)\in
\Real^n$. To show this it is enough to decompose system
(\ref{eq:intro:plant})--(\ref{eq:intro:CE_alg}) into the following
form
\begin{equation}\label{eq:intro:extended}
\begin{split}
\left(\begin{array}{l}
\dpsi\\
\dot{\tilde{\thetavec}}\\
\end{array}\right)
& = \left(\begin{array}{ll}
- \mathcal{\phi}_{\psi} & \mathcal{F}(\bfx,\thetavec,\hat{\thetavec})\\
\mathcal{A}_{\psi}(\psi,\bfx,\hat{\thetavec},t) & 0
\end{array}\right)
\left(\begin{array}{l}
\psi\\
\tilde{\thetavec}\\
\end{array}\right)\\
& + \left(\begin{array}{l}
1\\
0
\end{array}\right)\upsilon(t), \ \ \tilde{\thetavec}=\thetavec-{\hat{\thetavec}}
\end{split}
\end{equation}
where functions $\phi_{\psi}$,
$\mathcal{F}(\bfx,\thetavec,\hat{\thetavec})$,
$\mathcal{A}_{\psi}(\psi,\bfx,\hat{\thetavec},t)$ follow from
Hadamard lemma\footnote{In particular, these functions can be
calculated as follows
$\mathcal{F}(\bfx,\thetavec,\hat{\thetavec})=\int_0^{1}\frac{\pd
f(\bfx,\lambda \thetavec+(1-\lambda)\hat{\thetavec})}{\pd \lambda
\thetavec+(1-\lambda)\hat{\thetavec}} d\lambda$,
$\mathcal{\phi}_\psi=\int_0^1\frac{\pd \varphi(\lambda \psi)}{\pd
\lambda \psi} d\lambda$,
$\mathcal{A}_{\psi}(\psi,\bfx,\hat{\thetavec},t)=\int_0^{1}\frac{\pd
A(\lambda \psi,\bfx,\hat{\thetavec},t)}{\pd \lambda \psi}
d\lambda$}. Unlike in linear parametrization case, explicit
dependance of function
$\mathcal{F}(\bfx,\thetavec,\hat{\thetavec})$ on unknown
$\thetavec$ does not allow to compensate for uncertainty by
choosing appropriate function
$\mathcal{A}_{\psi}(\psi,\bfx,\hat{\thetavec},t)$ in
(\ref{eq:intro:extended}). Therefore, in general, it is necessary
to use additional control input $\upsilon(t)$ in order to ensure
that system (\ref{eq:intro:extended}) is Lyapunov stable and that
$\psi(\bfx)\rightarrow 0$ as $t\rightarrow \infty$. Success of
this strategy
is reported in \cite{Annaswamy99}, \cite{Lin},
\cite{Lin_2002_smooth}.


One step forward towards obtaining  non-dominating adaptive
control would be to reformulate the problem  as follows: design
function $\hat{\thetavec}(\bfx,t)$ such that
either
\begin{equation}\label{eq:intro:new_goal}
\lim_{t\rightarrow\infty}
f(\bfx,\thetavec)-f(\bfx,\hat\thetavec(\bfx,t))= 0,
\end{equation}
or
\begin{equation}\label{eq:intro:new_goal2}
f(\bfx,\thetavec)-f(\bfx,\hat\thetavec(\bfx,t))\in L_2
\end{equation}
hold\footnote{ Similar idea was proposed in \cite{Ortega02} as
adaptive ``root-searching" strategy}.
One possible way to achieve goal
(\ref{eq:intro:new_goal}), (\ref{eq:intro:new_goal2}) is to use
information about the difference
$f(\bfx,\thetavec)-f(\bfx,\hat\thetavec(\bfx,t))$ explicitly in
the adaptation algorithm rather than using external control
$\upsilon(t)$. For a class of nonlinear parameterized systems this
additional information can be introduced into adjustment schemes
 by mere structural changes in the adjustment law. In particular it is suggested in
\cite{t_fin_formsA&T2003}  to use adaptive algorithms in
differential-integral or {\it finite form} instead of using
adaptive algorithms in differential form. Extended system in this
case can be described as follows:
\begin{equation}\label{eq:intro:extended_2}
\left(\begin{array}{l}
\dpsi\\
\dot{\tilde{\thetavec}}\\
\end{array}\right)
= \left(\begin{array}{l}
- \mathcal{\phi}_{\psi} \ \ \ \mathcal{F}(\bfx,\thetavec,\hat{\thetavec})\\
0 \ - \mathcal{F}(\bfx,\thetavec,\hat{\thetavec})\alpha(\bfx,t)
\end{array}\right)
\left(\begin{array}{l}
\psi\\
\tilde{\thetavec}
\end{array}\right),
\end{equation}
where $\mathcal{F}(\bfx,\thetavec,\hat{\thetavec})\alpha(\bfx,t)$
is positive semi-definite time-varying matrix. Sufficient
conditions for existence of such algorithms are given in
\cite{ECC_2003},\cite{ALCOSP_2004}, where the problem is reduced
to solution of specific partial differential equation ( or {\it
explicit realizability condition}) by embedding the system into
the one of the higher order. While solution to this problem was
shown to exist for wide range of functions $\alpha(\bfx,t)$ and
$\psi(\bfx)$, it is not always possible to guarantee that
$\mathcal{F}(\bfx,\thetavec,\hat{\thetavec})\alpha(\bfx,t)$ is
positive semi-definite for  arbitrary $f(\bfx,\thetavec)$ and any
$\thetavec\in\Omega_\theta$, $\bfx(t_0)\in \Real^n$. As a result
of that Lyapunov stability of the whole system becomes
problematic, if not at all impossible, for general nonlinear
parameterizations.

These observations suggest that in both conventional problem
statement (\ref{eq:intro:error}), (\ref{eq:intro:CE_alg}) and
nonconventional ones (\ref{eq:intro:extended}),
(\ref{eq:intro:extended_2}) ensuring  Lyapunov stability for
adaptive control system in case of nonlinear parametrization
inevitably leads either to domination of the nonlinearity or to
additional restrictions on the class of nonlinear
parameterizations. In the other words, the problem of Lyapunov
stable and non-dominating adaptive control is ill-posed in
general. As a candidate for replacement of  Lyapunov stability one
could think of the {\it set attractivity}
\cite{Guckenheimer_2002,Milnor_85} of the goal manifold.
This concept allows us to design systems which being unstable in
Lyapunov sense have bounded solutions and also are capable of
reaching the goal asymptotically.
%
%
%
%
%
The main problem with this concept, however, is that
there is a number of conditions to check which critically depend
on precise knowledge of the vector fields of adaptive system. This
knowledge includes in particular properties of yet unknown
function $\hat{\thetavec}(\bfx,t,\delta)$.


Therefore perhaps the most reasonable and the least demanding
concept of the control goals  would be the notion of {\it
$\omega$-limit set}
\begin{defn} A point $p \in \Real^n$ is called an
$\omega$-limit point $\omega(\bfx(t,\bfx_0,t_0))$ of $\bfx_0\in
\Real^n$ if there exists sequence $\{t_i\}$,
$t_i\rightarrow\infty$, such that $\bfx(t,\bfx_0,t_0)\rightarrow
p$. The set of all limit points $\omega(\bfx(t,\bfx_0,t_0))$ is
the $\omega$-limit set of $\bfx_0$.
\end{defn}
Let, therefore, the control  goal be to ensure that for some
positive $\varepsilon$ the set
\begin{eqnarray}\label{eq:intro:goal_omega}
\Omega_\psi(\varepsilon)=\{\bfx\in \Real^n| \bfx:
|\psi(\bfx)|\leq\varepsilon\}
\end{eqnarray}
contains the {\it $\omega$-limit set} of $\Omega_\bfx$ for
non-autonomous system (\ref{eq:intro:plant}),
(\ref{eq:intro:control}) with $\upsilon(t)\equiv 0$. Hence, the
main question of our current study is the following: is there
exists a non-dominating adaptive scheme and reasonably large class
of parameterizations of the uncertainty such that all
$\Omega_\psi(\varepsilon)$ contains the $\omega$-limit set of the
adaptive system for any arbitrary small $\varepsilon>0$, all
trajectories of the system are bounded, and volume of the domain
of uncertainty is decreasing with time? The answer to this
question is provided in the next section.

\section{Main Results}

To begin with, let us  consider the case where function
$f(\cdot,\cdot)$ is parameterized by scalar $\theta\in
\Omega_{\theta}=[\underline{\theta},\overline{\theta}]\subset
\Real$, $\underline{\theta}<\overline{\theta}$. For each
$\theta\in \Omega_\theta$ and nonnegative $\Delta\in \Real_{\geq
0}$ we introduce the following equivalence relation
\[
\theta \sim \theta'\Leftrightarrow \ \
|f(\bfx,\theta)-f(\bfx,\theta')|\leq \Delta \ \forall \ \bfx\in
\Real^{n}
\]
and corresponding equivalence classes
$[{\theta}]_\Delta=\{\theta'\in \Omega_\theta|
\theta\sim\theta'\}$.

For the given functions $\varphi(\psi)$ and
$\alpha(t):\Real_{+}\rightarrow \Real$, $\alpha(t)\in C^1$ let us
define the following function
\[
S_{\delta}(\varphi(\psi),\alpha(t))=\left\{
                        \begin{array}{ll}
                         1, & |\varphi(\psi)+\alpha(t)|>\delta\\
                         0, & |\varphi(\psi)+\alpha(t)|\leq \delta
                        \end{array}
                 \right.
\]
With the function $S_\delta(\varphi(\psi(\bfx(t))),\alpha(t))$ we
associate the time sequence $\mathcal{T}=\underline{t}_0\leq
\overline{t}_0< \underline{t}_1<
\overline{t}_1<\dots<\underline{t}_i<\overline{t}_i<\underline{t}_{i+1}<\overline{t}_{i+1}<\dots$,
where
\begin{equation}\label{main:t_i}
\begin{split}
\underline{t}_0&= t_0  \\
\overline{t}_i&=\inf_{t\geq \underline{t}_i} \{t: \
|\varphi(\psi(\bfx(t)))+\alpha(t)|<\delta\} \\
\underline{t}_i&=\inf_{t\geq \overline{t}_{i-1}} \{t: \
|\varphi(\psi(\bfx(t)))+\alpha(t)|>\delta\}
\end{split}
\end{equation}
The elements of this sequence are time instances $\underline{t}_i$
(or $\overline{t}_i$) at which the sum
$\varphi(\psi(\bfx(t)))+\alpha(t)$ leaves (or enters) domain
$|\varphi(\psi(\bfx(t)))+\alpha(t)|\leq\delta$. We define that
$\underline{t}_0=\overline{t}_0$ if
$|\varphi(\psi(\bfx(t_0)))+\alpha(t_0)|<\delta$. Let us, in
addition, introduce function $\lambda$ with the following
properties:
\begin{equation}\label{main:lambda_single}
\begin{split}
&  \lambda:
\Real\rightarrow[\underline{\theta},\overline{\theta}], \
\lambda\in C^1, \
\im(\lambda)\supset[\underline{\theta},\overline{\theta}]\\
&  \forall \ s\in \Real, \theta\in\Omega_{\theta} \ \ \exists \
T, \tau(s)>0: \\
&  \theta=\lambda(s+\tau(s)), \ 0<\tau(s)<T
\end{split}
\end{equation}
An example of such function is
\[
\lambda(s)=\underline{\theta}+\overline{\theta}\frac{1}{2}(\sin(s)+1)
\]
As a candidate for $\hat{\thetavec}(\bfx,t,\delta)$ we choose the
following adaptation algorithm:
\begin{equation}\label{main:algorithm}
\begin{split}
\hat{\theta}(\bfx,t,\delta)&=\lambda(\hat{\theta}_0(\bfx,t,\delta)) \\
\hat{\theta}_0(\bfx,t,\delta)&=\gamma \left(\hat{\theta}_P(\bfx,t)+\theta_I(t)+C_{\theta}(t)\right) \\
\hat{\theta}_{P}(\bfx,t)&=
\psi(\bfx)\left(\alpha(t)+\frac{1}{2}\psi(\bfx)\right) \\
\dot{\hat{\theta}}_I&=S_{\delta}(\varphi(\psi),\alpha(t))(\psi(\bfx)\varphi(\psi)-\\
&  \psi(\bfx)(\xi_2 + b_1 \psi(\bfx)))\\
\alpha(t)&=(1, 0)(\xi_1, \xi_2)^{T}r\\
\left( \begin{array}{c}
 \dot{\xi}_1 \\
 \dot{\xi}_2
\end{array} \right) & = \left( \begin{array}{cc}
                0  & 1 \\
               a_1 & a_2
               \end{array} \right) \left( \begin{array}{c}
                                          \xi_1\\
                                          \xi_2
                                          \end{array} \right) +
                                          \left( \begin{array}{c}
                                                 b_1 \\
                                                 b_2
                                                 \end{array}
                                                 \right) \psi \\
&  b_1\neq0, \ a_1,a_2<0, \\
C_\theta(t)&=\left\{
            \begin{array}{c}
            \frac{1}{\gamma}\hat{\theta}(\overline{t}_{i-1})-\hat{\theta}_{I}(\overline{t}_{i-1})-\\
             \hat{\theta}_{P}(t), \ t\in (\overline{t}_{i-1},\underline{t}_{i})\\
           C_\theta(\overline{t}_{i-1}) + \hat{\theta}_{P}(\overline{t}_{i-1})- \\
            \hat{\theta}_{P}(\underline{t}_{i}), \ t\in
           [\underline{t}_{i},\overline{t}_i]
            \end{array}
            \right.
\end{split}
\end{equation}
Properties of algorithm (\ref{main:algorithm}) are summarized in
the following theorem:
\begin{thm}\label{theorem:nonlinear_single_dim} Let system (\ref{eq:intro:plant})
with control function (\ref{eq:intro:control}) and corresponding
error model (\ref{eq:intro:error}) be given. Let function
$\upsilon(t), \dot{\upsilon}(t)\in L_\infty$ and
$\|\upsilon(t)\|_\infty\leq \Delta$. Let in addition
function $f(\bfx,\thetavec)$ be bounded. Then\\
1) for any $\varepsilon>0$ and $\varphi\in
\mathcal{C}_\varphi(k)$, $k>0$ there exist functions
$\delta_0:\Real_+\rightarrow \Real_+, \ \delta_0\in C^0, \
\delta(0)=0$,
$\delta(\varepsilon,\Delta)=\delta_0(\varepsilon)+\Delta$, and
function $\hat{\theta}(\bfx,t,\delta)$, given by
(\ref{main:algorithm}) with arbitrary $\gamma\in \Real, \gamma>0$
and initial conditions
such that $\Omega_\psi(\varepsilon+\frac{\Delta}{k})$ contains the
$\omega$-limit set of system
(\ref{eq:intro:plant});\\
2) all trajectories of the system are bounded and solutions
$\bfx(t,\bfx_0,t_0)$ converge into the domain
$\Omega_\psi(\varepsilon+\frac{\Delta}{k})$ in finite time;\\
3) if for any $\theta\in\Omega_\theta$ there exist constants
$T_1>0$, $M>2\delta_0(\varepsilon)+\Delta>0$ and function
$\tau(t):\Real \rightarrow (0,T_1)$ such that
\begin{equation}\label{nonlinear_PE}
\begin{split}
& |f(\bfx(t+\tau(t)),\theta)-f(\bfx(t+\tau(t)),\hat{\theta})|>M \\
&  \forall \ \hat{\theta}\in \  \Omega_{\theta}\ \backslash \
\mathcal{U}_\epsilon([\theta])
\end{split}
\end{equation}
then $\hat{\theta}$ converges into
$\mathcal{U}_\epsilon([\theta])$ in finite time.\\
{\it Proof of the Theorems are given in Appendix.}
\end{thm}
 Theorem
\ref{theorem:nonlinear_single_dim} states that for arbitrary
$C^1$-smooth and bounded function $f(\bfx,\theta)$ there exists
non-dominating adaptation algorithm in class
$\mathcal{C}_\varphi(k)$. To see this it is enough to let
$\Delta=0$. In presence of the unknown perturbation $\upsilon(t)$
we guarantee convergence of the trajectories $\bfx(t,\bfx_0,t_0)$
to arbitrary small neighborhood of $\Omega_\psi(\frac{\Delta}{k})$
subject to the choice of $\delta_0>0$.

Adaptation algorithm (\ref{main:algorithm}) ensures boundedness of
the solutions in the extended system and furthermore, under
assumption (\ref{nonlinear_PE}), it guarantees convergence of the
estimates $\hat{\theta}$ to arbitrary small neighborhood
$\mathcal{U}_\epsilon([\theta]_\Delta)$ of the equivalence class
(set) $[\theta]_\Delta$. In case $[\theta]_\Delta=\{\theta\}$ it
guarantees convergence to the small neighborhood of the actual
value of $\theta$.

Condition (\ref{nonlinear_PE}), which we require for convergence
of the parameter $\hat{\theta}$ into $\mathcal{U}_\epsilon$, can
be regarded as a new version of {\it nonlinear persistent
excitation} \cite{Cao_2003}. Our condition, however, is more easy
to verify. In addition, this condition is consistent with {\it
linear persistent excitation} condition \cite{Narendra89}:
$\exists \ T>0, \ \rho>0: \ \int_{t}^{t+T}\bfx(\tau)\bfx(\tau)^T
d\tau
> \rho I_n$. Indeed $\int_{t}^{t+T}\bfx(\tau)\bfx(\tau)^T d\tau
> \rho I_n \Rightarrow  \forall \thetavec\neq\thetavec' \
(\thetavec-\thetavec')^{T}\times$
$\int_{t}^{t+T}\bfx(\tau)\bfx(\tau)^T d\tau
(\thetavec-\thetavec')>\rho \|\thetavec-\thetavec'\|^2\Rightarrow$
$|(\thetavec-\thetavec')^T\bfx(t_1)|>\frac{\rho}{T}\|\thetavec-\thetavec'\|,
\ t_1\in[t,t+T] \Rightarrow$
$|(\thetavec-\thetavec')^T\bfx(t_1)|>M, \ \forall
\theta\in\Omega_\theta\backslash\mathcal{U}_{\epsilon}, \
M=\frac{\rho\epsilon}{T}$

Let us generalize the statements of Theorem
\ref{theorem:nonlinear_single_dim} to the case where $\thetavec\in
\Omega_\theta\subset \Real^d$. For that reason we introduce the
following assumption:
\begin{assume}\label{assume:multi-dim} Let $\Omega_\theta$  be bounded and there exist
$C^{1}$-smooth function $\etavec:
[\underline{\theta},\overline{\theta}]\rightarrow \Real^{d}$ such
that that for any $\thetavec\in\Omega_\theta$ there exists
$\lambda^{\ast}(\thetavec)\in
[\underline{\theta},\overline{\theta}]$:
$|f(\bfx,\thetavec)-f(\bfx,\etavec(\lambda^{\ast})|\leq \Delta, \
\forall \ \bfx\in \Real^{n}$
\end{assume}
Applicability of algorithms (\ref{main:algorithm}) to
multi-dimensional $\thetavec$ then follows explicitly from Theorem
\ref{theorem:nonlinear_single_dim}.
\begin{thm}\label{theorem:nonlinear_multi_dim}  Let system (\ref{eq:intro:plant})
with control function (\ref{eq:intro:control}) and corresponding
error model (\ref{eq:intro:error}) be given, function
$f(\bfx,\thetavec)$ be bounded, and Assumption
\ref{assume:multi-dim} hold. Then statements 1) -- 3) of Theorem
\ref{theorem:nonlinear_single_dim} follow.
\end{thm}

Adaptation algorithm (\ref{main:algorithm}) can be considered as a
nonlinear dynamical system which, however, is not internally
globally stable in Lyapunov sense. The algorithm, nonetheless,
does not lead to unbounded growth of its internal state, neither
it fails to ensure reaching of the control goal for arbitrary
initial conditions $\bfx(t_0)$, $\hat{\theta}_{I}(t_0)$,
$\xivec(t_0)$. The properties of this algorithm essentially rely
on two ideas: {\it monotonic evolution} of $\hat{\theta}_0(t)$,
and {\it multiple equilibriums} in the corresponding differential
equation for ${\hat{\theta}}_0$ given by (\ref{main:d_theta_eq})
in Appendix.

Multiple equilibriums of  are guaranteed by function
$\lambda(\cdot)$ defined by (\ref{main:lambda_single}) and
invariance of $\hat{\theta}_0(t)$ on the following set $
\{\bfx,\hat{\theta}_0| \ \bfx,\hat{\theta}_0: \
f(\bfx,\theta)-f(\bfx,\lambda(\hat{\theta}_0))=0 \}$. Monotonic
character of  ${\hat{\theta}}_0(t)$ and multiple equilibriums
ensure existence of the limit $\lim_{t\rightarrow
\infty}{\hat{\theta}}_0(t)=\hat{\theta}_{0,\infty}$  for any
initial conditions $\bfx(t_0)$, $\hat{\theta}_{I}(t_0)$,
$\xivec(t_0)$. This fact is used in the proof of Theorem
\ref{theorem:nonlinear_single_dim} to show convergence of the
trajectories $\bfx(t,\bfx_0,t_0)$ to the set specified by
(\ref{eq:intro:goal_omega}). Notice also that existence of these
multiple equilibriums  in case of persistent excitation for
$\upsilon(t)\equiv 0$ ensure that  the set
$\mathcal{U}_\epsilon([\theta])$ becomes globally attractive
w.r.t. function $\hat{\theta}(t,\bfx,\delta)$.

\begin{figure}\label{fig:multistable}
\begin{center}
\includegraphics[width=\columnwidth]{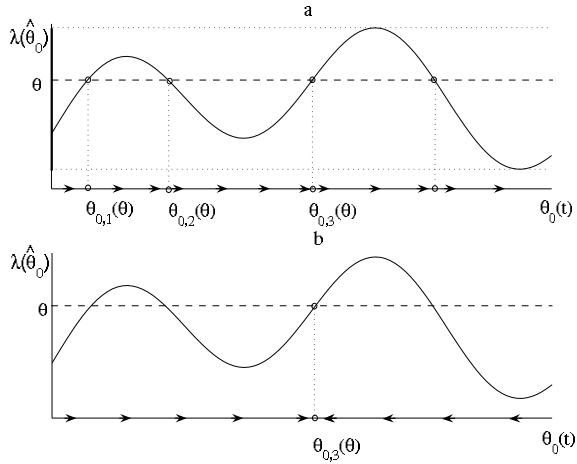}
\end{center}
\begin{center}
\caption{Adaptation with multiple equivalent equilibriums (plot a)
vs. adaptation with single asymptotically stable equilibrium (plot
b)}
\end{center}
\vskip 65mm
\end{figure}
\begin{figure}\label{fig:trajectory}
\begin{center}
\includegraphics[width=80pt]{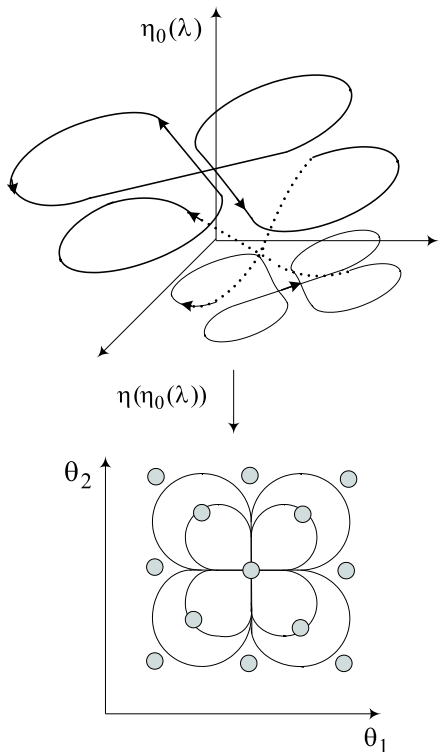}
\end{center}
\vskip 7mm
\begin{center}
\caption{Trajectory $\etavec(\lambda)$ (the bottom panel) as
projection of the smooth curve. Gray circles form the gird in the
parameters space which induces the curves satisfying Assumption
\ref{assume:multi-dim}.}
\end{center}
\end{figure}
Conceptual difference between our method and conventional
Lyapunov-based design is illustrated with Fig. 1. In Fig. 1 the
upper plot depicts the solution curve
$\hat{\theta}_0(t,\hat{\theta}_0(t_0),t_0)$ of differential
equation (\ref{main:d_theta_eq}), were the arrows point towards
increase of the independent variable along the curve. For the
given value of $\theta$ and initial condition
$\hat{\theta}_0(t_0)$ function $\lambda(\hat{\theta_0})$ generates
infinitely many equilibria $\hat{\theta}_{0,i}$, $i\in \Natural$.
If the perturbation is applied to the system the system will
eventually escape its current equilibrium (for instance,
$\hat{\theta}_0=\hat{\theta}_{0,1}$) equilibrium and move along
the axis $\hat{\theta}_0$. Due to the monotonicity of
$\hat{\theta}_0(t,\hat{\theta}_0(t_0),t)$ with respect to $t$ it
eventually reaches a neighborhood of the point
$\hat{\theta}_0=\hat{\theta}_{0,2}$ and stops there if the
perturbation is released. In order to prevent unbounded growth of
$\hat{\theta}_0$ under persistent perturbations it may be
necessary to change the sign of $\gamma$ in
(\ref{main:d_theta_eq}) upon solution $\hat{\theta}_0(t)$ reaches
certain bounds.

In the lower plot we have shown the typical would be solution
curves of Lyapunov asymptotic stable estimator. Disturbances will
move the solution away from the equilibrium and the system will
restore its original state after the disturbing terms vanish.
Despite obvious advantages of this type of behavior, i.e. small
perturbations induce small deviations, in adaptive control
problems the parametric perturbations are usually large.
Furthermore, as it has been already pointed out, gentle control
ensuring Lyapunov stability for arbitrary nonlinear parameterized
system is problematic.

One of the distinctive features of algorithms
(\ref{main:algorithm}) is that they can in principle take into
account information about (available a-priori) distribution of the
unknown quasi-stationary $\thetavec$ as a function of time. This
information is then to be accounted for by the choice of functions
$\etavec(\cdot)$ and $\lambda(\cdot)$. Illustration to this choice
is given in Figure 2. In Fig. 2 the curve $\etavec(\lambda)$ is
designed to visit neighborhood of the center more frequently ($16$
times per period) than other points (only $2$ times per period) of
the domain. Despite the problem of choosing these curves
$\etavec(\lambda)$ (which fit given distributions of $\thetavec$)
is not trivial, such optimal choice, if successful, can provide
room for further enhancements of performance in the adaptive
systems. The second step would be to adjust or tune functions
$\etavec(\lambda)$ adaptively thus enabling self-tuning of the
adaptive algorithm itself. These topics, however, are beyond the
goals of our current study.





\section{Conclusion}

In this paper we have proposed new technique for adaptive control
of nonlinear dynamical systems with nonlinear parametrization. In
contrast to conventional concept of  Lyapunov stable  adaptive
control, and as a result domination of the nonlinearity by
high-gain feedbacks, we use adaptation schemes which are not
stable in Lyaponov sense. Yet, these algorithms guarantee reaching
of arbitrary small neighborhood of the desired target set.
Moreover the resulting control function is non-dominating.

The ideology beyond our method is somewhat similar to the one
introduced in \cite{Ilchman_97,Pomet92}. The results, however, are
substantially different. First, we do not require neither
exponential stability nor asymptotic stability  of the target
dynamics. Second, the speed of adaptation in our case is not to be
slowed down with time. Last but not least is that the method can
be used to identify nonlinear systems of rather general class
without requesting for linearization of the nonlinearities.

Technique for adaptive control introduced in our paper rather well
coincides with the thesis of "nonlinear philosophy for nonlinear
systems" declared in \cite{Fradkov_2000}. Furthermore, it fits
well empirical observations that complex natural systems poses
clearly detectable patterns of instability in the empirical data.

Despite we proved only that the target set will be reached in the
system with our algorithms, we hope that some type of robustness
w.r.t. unknown external perturbations can be also shown. In fact,
robustness of the system to unmodeled dynamics with known
$L_\infty$ norm can be easily ensured by enlarging the value of
$\delta$ in our algorithms. Whether or not robust behavior can be
achieved by choice of another parameters like functions
$\etavec(\cdot)$ or $\lambda(\cdot)$ will be the topics of our
future study.

\bibliography{general_parametrization_conf}

\section{Appendix}

{\it Proof of Theorem \ref{theorem:nonlinear_single_dim}.} Let us
consider system (\ref{main:algorithm}) and calculate derivative
$\dot{\hat{\theta}}_0$ with respect to independent variable  $t$:
\begin{equation}\label{main:d_theta_eq}
\begin{split}
\dot{\hat{\theta}}_0&=\gamma
S_\delta(\psi,\alpha(t))(\varphi(\psi)+\dpsi)(\varphi(\psi)+\alpha(t)),\\
& \ \hat{\theta}_0(t_0)\in\Real
\end{split}
\end{equation}
To proceed further we will need auxiliary Lemma
\ref{boundedness_derivative} and Lemma \ref{derivative}.

\begin{lem}\label{boundedness_derivative} Consider system
(\ref{eq:intro:plant}), (\ref{eq:intro:error}). Let function
$f(\bfx,\thetavec)$ in (\ref{eq:intro:error}) be bounded w.r.t.
$\bfx,\thetavec$ and $\upsilon(t), \dot{\upsilon}(t)\in L_\infty$.
Then $\hat{\thetavec}(t),\dot{\hat\thetavec}(t)\in L_\infty$ imply
that
 $\ddot{\psi}\in L_\infty$.
\end{lem}
{\it Proof of Lemma \ref{boundedness_derivative}}. First, we
observe that $\psi(\bfx)\in L_\infty$. This follows immediately
from boundedness of $f(\bfx,\thetavec)$ and the fact that function
$\varphi$ belongs to the class $\mathcal{C}_\varphi$ specified by
(\ref{eq:intro:class_phi}). Therefore, according to assumption
(\ref{eq:intro:psi_assume}) on function $\psi(\bfx)$, state
$\bfx\in L_\infty$. Moreover, given that $\psi\in L_\infty$ and
$\varphi\in C^1$, we can derive from (\ref{eq:intro:error}) that
$\dpsi\in L_\infty$. Let us consider $\ddot{\psi}$:
$\ddot{\psi}=\left(\frac{\pd f(\bfx,\thetavec)}{\pd \bfx}-
\frac{\pd f(\bfx,\hat{\thetavec})}{\pd
\bfx}\right)(\bff(\bfx,\thetavec)+\bfg(\bfx)u)+\frac{\pd
f(\bfx,\hat\thetavec)}{\pd \hat{\thetavec}}\dot{\hat{\thetavec}} +
\frac{\pd \varphi(\psi)}{\pd \psi}\dpsi+\dot{\upsilon}(t)$.
Derivatives $\frac{\pd f(\bfx,\thetavec)}{\pd \bfx}$, $\frac{\pd
f(\bfx,\hat{\thetavec})}{\pd \bfx}$ are continuous functions with
respect to $\bfx$ as  $\psi\in C^2, \bff(\bfx,\thetavec)\in C^1$.
Therefore they are bounded as $\bfx, \hat{\thetavec}, \thetavec
\in L_\infty$. In addition, continuity of $\frac{\pd
f(\bfx,\hat\thetavec)}{\pd \hat{\thetavec}}$ and $\frac{\pd
\varphi}{\pd \psi}$ results in their boundedness. To conclude the
proof it is enough to notice that function $u$ is bounded for
bounded $\bfx,\hat{\thetavec}$. {\it The lemma is proven.}

\begin{lem}\label{derivative}
Let function $\psi(t)$ be given and its second time-derivative is
bounded: $|\ddot{\psi}|<\beta_1$, $\beta_1>0$. Then there exists
differential filter
\begin{equation}\label{lf1}
\begin{split}
\left( \begin{array}{c}
 \dot{\xi}_1 \\
 \dot{\xi}_2
\end{array} \right) & =  \left( \begin{array}{cc}
                0  & 1 \\
               a_1 & a_2
               \end{array} \right) \left( \begin{array}{c}
                                          \xi_1\\
                                          \xi_2
                                          \end{array} \right) +
                                          \left( \begin{array}{c}
                                                 b_1 \\
                                                 b_2
                                                 \end{array}
                                                 \right) \psi
                                                 \\
y &= \left( c_1, \ 0 \right)  \left( \xi_1, \ \xi_2\right)^T,
\end{split}
\end{equation}
where $b_1\neq0, a_1,a_2<0$, and time $t_1>0$ such that for any
positive constant $\epsilon>0$ the following estimate holds:
$|y(t)-\dot{\psi}(t)|\leq \epsilon \ \forall \ t>t_1$. In
particular, if $c_1b_1=- a_1$, $b_2=a_2 b_1$ then filter output
satisfies the following inequality:
\begin{equation}\label{filter_error}
|y(t)-\dot{\psi}(t)|\leq |\frac{a_2\beta_1}{a_1}|+\delta(t),
\end{equation}
where $\delta(t)$ decays exponentially fast to the origin.
\end{lem}

{\it Proof of Lemma \ref{derivative}}. Let $b_2=a_2 b_1$. Then
system (\ref{lf1}) has the following equivalent representation:
\begin{equation}\label{li1}
\begin{split}
\left( \begin{array}{c}
 \dot{\xi}_1 \\
 \dot{\xi}_2
\end{array} \right) & =  \left( \begin{array}{cc}
               0 & 1 \\
               a_1 & a_2
               \end{array} \right) \left( \begin{array}{c}
                                          \xi_1\\
                                          \xi_2
                                          \end{array} \right) +
                                          \left( \begin{array}{c}
                                                 0 \\
                                                 b_1
                                                 \end{array}
                                                 \right) \dot{\psi}
                                                 \\
y &= \left( c_1 \ 0 \right) \left(\xi_1, \ \xi_2\right)^T
\end{split}
\end{equation}
Denoting
\[
A=\left( \begin{array}{cc}
               0 & 1 \\
               a_1 & a_2
               \end{array} \right), \ \ \mathbf{b}_1=\left( \begin{array}{c}
                                                 0 \\
                                                 b_1
                                                 \end{array}
                                                 \right), \ \
                                                 \mathbf{c}=\left( \begin{array}{c}
                                                 c_1 \\
                                                 0
                                                 \end{array}
                                                 \right),
\]
we may write the filter equations in a vector-matrix form:
\begin{equation}\label{li2}
\begin{split}
\dot{\xivec}&=A\mathbf{\xivec}+\mathbf{b}_1\dot{\psi} \\
 y&=\mathbf{c}^T\xivec
 \end{split}
\end{equation}
Consider output $y$ of system (\ref{li2}):
\begin{equation}\label{l3}
 y=\mathbf{c}^T\left( e^{At}\mathbf{\xivec}_0+e^{At}\int_0^t e^{-A \tau}\mathbf{b}_1 \dot{\psi}(\tau) d\tau
 \right).
\end{equation}
Matrix $A$ is Hurwitz and,  therefore, invertible. Hence taking
into account existence of $\ddot{\psi}(t)$ and that matrices
$A^{-1}$, $e^{At}$ commute, we can rewrite
 equality  (\ref{l3}) as follows: $y = \mathbf{c}^T( e^{At}\mathbf{\xivec}_0-
 e^{At}A^{-1}e^{-At}\mathbf{b}_1\dot{\psi}(t) {|}_0^t   +
e^{At}A^{-1}\int_0^t e^{-A\tau}\mathbf{b}_1 \ddot{\psi}(\tau)
 d\tau)$. Hence $y = \mathbf{c}^T( e^{At}\mathbf{\xivec}_0-A^{-1}
 \mathbf{b}_1\dot{\psi}(t) + A^{-1}e^{At}\mathbf{b}_1\dot{\psi}(0)+A^{-1}\int_0^te^{A(t-\tau)}\mathbf{b}_1 \ddot{\psi}(\tau)
 d\tau )$.
Consider the following difference $|y(t)-\dot{\psi}(t)|$:
\begin{equation}
\begin{split}
& |y(t)-\dot{\psi}(t)|=|\mathbf{c}^T(
e^{At}\mathbf{\xivec}_0-A^{-1}
 \mathbf{b}_1\dot{\psi}(t) +
A^{-1}\times\\
& e^{At}\mathbf{b}_1\dot{\psi}(0)+ A^{-1}\int_0^t
e^{A(t-\tau)}\mathbf{b}_1 \ddot{\psi}(\tau)
 d\tau ) - \dot{\psi}(t)|  \\
&  \leq |\mathbf{c}^T ( e^{At}\xivec_0-A^{-1}
 \mathbf{b}_1\dot{\psi}(t) + A^{-1}e^{At}\mathbf{b}_1\dot{\psi}(0)
 )\\
 & -  \dot{\psi}(t)| + | \mathbf{c}^T A^{-1}\int_0^te^{A(t-\tau)}\mathbf{b}_1\beta_1 d\tau|  \\
 &  =  |\mathbf{c}^T( e^{At}\xivec_0-A^{-1}
 \mathbf{b}_1\dot{\psi}(t) + A^{-1}e^{At}\mathbf{b}_1\dot{\psi}(0)
 )\\
 & - \dot{\psi}(t)| + |\mathbf{c}^T A^{-2}e^{At}\mathbf{b}_1\beta_1-\mathbf{c}^T
A^{-2}\mathbf{b}_1\beta_1|,\nonumber
\end{split}
\end{equation}
where $A^{-2}=A^{-1}A^{-1}$. Let us denote
$\delta(t)=|\mathbf{c}^T e^{At}\xivec_0|+|\mathbf{c}^T
A^{-1}e^{At}\mathbf{b}_1\dot{\psi}(0)|+|\mathbf{c}^T
A^{-2}e^{At}\beta_1|$. Term $\delta(t)\rightarrow0$ at
$t\rightarrow\infty$ as matrix $A$ is Hurwitz. Therefore for any
$\delta_1>0$ there exists time $t_1>0$ such that for any $t>t_1$
the following inequality holds:
\begin{equation}\label{l6}
\begin{split}
  |y(t)-\dot{\psi}(t)|&\leq |(\mathbf{c}^T A^{-1}
 \mathbf{b}_1 +1 )\dot{\psi}(t)|+\\
 &|(\mathbf{c}^TA^{-2}\mathbf{b}_1\beta_1|+ \delta_1.
 \end{split}
\end{equation}
Let us consider term $\mathbf{c}^T A^{-1}
 \mathbf{b}_1$ in (\ref{l6}). Matrix
 \[
A^{-1}=\left( \begin{array}{cc}
         0 & 1 \\
        a_1 & a_2
        \end{array}
        \right)^{-1}=\frac{1}{a_1}\left( \begin{array}{cc}
                             -a_2& 1\\
                              a_1&0
                            \end{array} \right)
 \]
 Then
 \[
\mathbf{c}^T A^{-1}\mathbf{b}_1=(c_1 \ 0) A^{-1} \left(
\begin{array}{c}
                                       0 \\
                                       b_1
                                       \end{array}
                                       \right)=\frac{b_1 \ c_1}{a_1}
 \]
Notice that $\frac{b_1 \ c_1}{a_1}=-1$ due to the lemma
assumptions. Therefore, using inequality (\ref{l6}) we can derive
the following estimate:
\[
|y(t)-\dot{\psi}(t)|\leq |(\mathbf{c}^TA^{-2}\mathbf{b}_1\beta_1|+
\delta_1,
\]
where
\[
A^{-2}=\frac{1}{a_1^2}\left( \begin{array}{cc}
                      a_2^2+a_1& -a_2 \\
                      -a_2 a_1 & 0
                      \end{array} \right).
\]
Hence
\begin{equation}\label{l7}
\begin{split}
 & |y(t)-\dot{\psi}(t)|\leq \left|(c_1 \ 0) A^{-2} \left(
\begin{array}{c}
                                                 0\\
                                                 b_1
                                                 \end{array}
                                \right)\beta_1 \right| + \delta_1\\
& =|\frac{a_2\beta_1}{a_1}|+\delta_1.
\end{split}
\end{equation}
Inequality (\ref{l7}) proves  Lemma \ref{derivative}.

Lemmas \ref{boundedness_derivative},\ref{derivative} will allow us
to show that $\hat{\theta}(\bfx,t,\delta)$,
$\delta=\delta_0+\Delta$, $\delta_0\in \Real_+$ has a limit as
$t\rightarrow\infty$. First we observe that
$\hat{\thetavec}(\bfx,t,\delta)$ is bounded. Moreover
$\hat{\thetavec}(\bfx,t,\delta)$ is continuous and has bounded
time-derivative. Therefore, according to Lemma
\ref{boundedness_derivative} derivative $\ddot{\psi}$ is bounded.
Therefore, as it follows from Lemma \ref{derivative}, there exist
parameters $a_1,a_2,b_1,b_2$ of filter (\ref{li1}) such that
$|\alpha(t)-\dpsi(t)|<\delta_0/4 + |\delta_1(t)|$, where
$\delta_1(t)\rightarrow 0$ as $t\rightarrow\infty$. In particular,
there exist time $t_1>0$ such that $|\delta_1(t)|<\delta_0/4$ for
any $t>t_1$. Notice, that function
${\hat{\theta}}_0(\bfx,t,\delta)$ is bounded for $t\leq t_1$ (as a
sum and integral of bounded functions in time). Let us show that
for any $t>t_1$ function $\hat{\theta}_0(\bfx,t,\delta)$ is
monotonically increasing with respect to time $t$.

Let us consider function $\alpha(t)$, $t>t_1$. Taking into account
that $|\alpha(t)-\dpsi(t)|<\delta_0/2$ for any $t>t_1$ we can
estimate function $\alpha(t)$ in the following way:
\begin{equation}\label{main:alpha_vaphi_est}
\alpha(t)=\dpsi(t)+\mu(t), \  |\mu(t)|\leq \delta_0/2
\end{equation}
Therefore, for any $\alpha(t),\psi$:
$|(\alpha(t)+\varphi(\psi))|>\delta$ we have
\begin{equation}\label{main:d_theta_0_mon}
\begin{split}
&
(\dpsi+\varphi(\psi))(\alpha(t)+\varphi(\psi))=(\alpha+\varphi(\psi)-
\\
&  \mu(t))(\alpha(t)+\varphi(\psi))=
(\alpha(t)+\varphi(\psi))^2-\\
& \mu(t)(\alpha(t)+\varphi(\psi))\geq
|(\alpha(t)+\varphi(\psi))|\delta
-\\
&\frac{\delta_0}{2}|(\alpha(t)+\varphi(\psi))|\geq
\frac{\delta_0^2}{2}>0
\end{split}
\end{equation}
Hence, according to equations (\ref{main:d_theta_eq}),
(\ref{main:d_theta_0_mon}) we can conclude that function
$\hat{\theta}_0(t,\bfx,\delta)$ is monotonic, continuous and not
decreasing with time.

Let us show that, function $\hat{\theta}_0(\bfx,t,\delta)$ is also
bounded from above. Assume that for any $D_{\theta}>0$ there
exists $t_2>0$ such that
\begin{equation}\label{main:hypothesis_unbound}
|\hat{\theta}_0(\bfx(t_1),t_1,\delta)-\hat{\theta}_0(\bfx(t_2),t_2,\delta)|>D_\theta
\end{equation}
Let us select $D_\theta>T$, where the values of  $T$ are specified
in (\ref{main:lambda_single}). Function
$\hat{\theta}_0(\bfx(\tau),\tau,\delta)$ is continuous w.r.t.
$\tau$ and bounded for any $\tau\in[t_1,t_2]$(boundedness of
$\hat{\theta}_0(\bfx(\tau),\tau,\delta)$ follows from boundedness
of its derivative over finite time interval). Therefore, according
to the intermediate value theorem,  for any  $0< T_0< D_\theta$
there exists $t_3\in[t_1,t_2]$ such that
$\hat{\theta}_0(\bfx(t_3),t_3,\delta)=\hat{\theta}_0(\bfx(t_1),t_1,\delta)+T_0$.
Taking into account (\ref{main:lambda_single}) we can conclude
that for any $\theta\in[\underline{\theta},\overline{\theta}]$
there exists $t_3>0$ such that
$\theta=\lambda(\hat{\theta}_0(\bfx(t_3),t_3,\delta)), \
t_3\in[t_1,t_2]$.  According to the properties of function
$\hat{\theta}_0(\bfx,t,\delta)$, however, derivative
$\dot{\hat{\theta}}_0(t)\equiv 0$ for any $t\geq t_3$ as
\begin{equation}
\begin{split}
&
|\alpha(t_3)+\varphi(\psi(\bfx(t_3)))|=|\dot{\psi}+\varphi(\psi(\bfx(t_3)))+\mu(t_3)|\nonumber
\\
&
=|f(\bfx,\lambda(\hat{\theta}_0(\bfx(t_3),t_3,\delta))-f(\bfx,\theta))+\upsilon(t_3)+\nonumber\\
& \mu(t_3)|\leq |\upsilon(t)|+|\mu(t)|\leq
\delta_0/2+\Delta<\delta \ \forall \ t\geq t_3 \nonumber
\end{split}
\end{equation}
Therefore
$\hat{\theta}_0(\bfx(t_2),t_2,\delta)=\hat{\theta}_0(\bfx(t_3),t_3,\delta)\leq
\hat{\theta}_0(\bfx(t_1),t_1,\delta)+T$ which contradicts to
(\ref{main:hypothesis_unbound}).

So far we have shown that $\hat{\theta}_0(\bfx(t),t,\delta)$ is
continuous, monotonic and bounded function for any $t>t_1$.
Therefore, applying Bolzano-Weierstrass,  theorem we can  derive
that there exists the following limit
\begin{equation}\label{main:theta_0_limit}
\lim_{\tau\rightarrow\infty}\hat{\theta}_{0}(\bfx(\tau),\tau,\delta)=\hat{\theta}_{0,\infty}
\end{equation}
Taking into account equality (\ref{main:d_theta_eq}) we can
rewrite function $\hat{\theta}_0(\bfx(t),t,\delta)$ for $t>t_1$ as
follows:
\begin{equation}\label{main:theta_0_int}
\begin{split}
&
\hat{\theta}_0(\bfx(t),t,\delta)=\hat{\theta}_0(\bfx(t_1),t_1,\delta)+\\
&
\int_{t_1}^{t}S_\delta(\varphi(\psi(\tau)),\alpha(\tau))(\dpsi(\tau)+\varphi(\psi(\tau)))\times\\
&(\varphi(\psi(\tau))+\alpha(\tau))d\tau
\end{split}
\end{equation}
Equations (\ref{main:theta_0_int}), (\ref{main:theta_0_limit}) and
inequality (\ref{main:d_theta_0_mon}) result in the following
estimate $0<\int_{t_1}^{\infty}S_\delta(\varphi(\psi(\tau)),$
$\alpha(\tau))\frac{\delta^{2}}{2}d\tau=
\sum_{i=0}^{\infty}(\overline{t}_i-\underline{t}_i)\frac{\delta^2}{2}<\infty$.
The last inequality implies that
$\lim_{i\rightarrow\infty}\overline{t}_i-\underline{t}_i=0$. Let
us consider $|\varphi(\psi(t))+\alpha(t)|$ for $t\in
[\underline{t}_i,\overline{t}_i]$:
\[
|\varphi(\psi(t))+\alpha(t)|\leq
|\varphi(\psi(\underline{t}_i))+\alpha(\underline{t}_i)|+(\overline{t}_i-\underline{t}_i)D_{|\cdot|},
\]
where
$D_{|\cdot|}=\max_{t\in[\underline{t}_i,\overline{t}_i]}\frac{d}{dt}|\varphi(\psi(t))+\alpha(t)|$.
Function $D_{|\cdot|}$ is the bounded function of time and
therefore
\begin{equation}\label{main:psi_limit}
\limsup_{t\rightarrow\infty}|\varphi(\psi(t))+\alpha(t)| = \delta
\end{equation}
Taking (\ref{main:psi_limit}) into account we can derive that
dynamics of function $\psi(t)$ for $t>t_1$ can be described as
follows:
\begin{equation}\label{main:dpsi_t_4}
\begin{split}
\dpsi&=\mu(t)-\varphi(\psi), \
\limsup_{t\rightarrow\infty}|\mu(t)|\leq \frac{1}{2}\delta_0 +\\
& \delta = \frac{3}{2}\delta_0 + \Delta
\end{split}
\end{equation}
According to the definition of $\limsup$ there exists a time
instant $t_4$ such that $|\mu(t)|<2\delta_0+\Delta$ for any
$t>t_4$. Let us consider the following function
\begin{equation}
\begin{split}
V(\psi,\nu)&=\int_{0}^{\psi(t)}\phi(\xi,\nu)d\xi, \\
\phi(\xi,\nu)&=\left\{
                 \begin{array}{cc}
                 \xi-\nu, & \xi>\nu \\
                 0, & -\nu \leq \xi\leq \nu\\
                 \xi+\nu, & \xi<-\nu
                 \end{array}
                 \right., \ \nu\geq 0 \nonumber
\end{split}
\end{equation}
where $\psi(t)$ is the solution of system (\ref{main:dpsi_t_4})
for $t>t_4$  with initial condition $\psi(t_4)$. Function
$\varphi(\psi)\in \mathcal{C}_\varphi(k)$, therefore
$\varphi(\psi)\psi\geq k \psi^{2}$. Let us select
$\nu=(2\delta_0+\Delta)/k$ and consider derivative $\dot{V}$:
\begin{equation}
\begin{split}
\dot{V}&=\phi(\psi,\nu)\dpsi=\phi(\psi,\nu)(\mu(t)-\varphi(\psi))\\
& \leq \phi(\psi,\nu)(\mu(t)-k\psi)\leq 0.\nonumber
\end{split}
\end{equation}
Moreover, it follows from Barbalat's lemma that
$\phi(\psi,\nu)\rightarrow 0$ as $t\rightarrow\infty$. This
automatically implies that there exist time instant $t_5$ such
that $\|{\psi(\bfx(t))}\|_\infty\leq (3\delta_0+\Delta)/k$ for any
$t>t_5$.

So far we have shown that for any $\delta>\Delta+\delta_0$,
$\delta_0>0$, $\varphi\in \mathcal{C}_\varphi(k)$ and
$\bfx_0\in\Omega_\bfx$ there exists algorithm
(\ref{main:algorithm}) with arbitrary $\gamma>0$  and initial
conditions such that $\Omega_\psi(2\delta_0/k+\Delta/k)$ contains
the $\omega$-limit set of system (\ref{eq:intro:plant}),
(\ref{eq:intro:control}) with $\upsilon(t): \
\|\upsilon(t)\|_\infty\leq \Delta$. Furthermore, we have shown
that trajectories $\bfx(t,\bfx_0,t_0)$ converge into the domain
$\Omega_\psi(3\delta_0/k+\Delta/k)$ in finite time. Therefore in
order to complete the proof, given arbitrary $\varepsilon>0$, it
is enough to pick $\delta_0<1/3 k\varepsilon$. Hence, statements 1
and 2 of the theorem are proven.

Let us show that statement 3 of the theorem holds. Notice that, in
order to complete the proof it is enough to show that
$\lambda(\hat{\theta}_{0,\infty})\in
\mathcal{U}_\epsilon([\theta]_\Delta)$, where
$\hat{\theta}_{0,\infty}$ is defined as in
(\ref{main:theta_0_limit}). Let us assume that
$\lambda(\hat{\theta}_0(\bfx,t_1,\delta))$ does not belong to
$\mathcal{U}_\epsilon([\theta]_\Delta)$, and consider equality
(\ref{main:alpha_vaphi_est}) for $t>t_1$. According to
(\ref{main:alpha_vaphi_est}) the following equality holds:
\begin{equation}
\begin{split}
&\dpsi+\varphi(\psi)+\mu(t)=\alpha(t)+\varphi(\psi)=f(\bfx,\theta)-f(\bfx,\hat{\theta})\\
& + \upsilon(t)+\mu(t), \ |\mu(t)+\upsilon|\leq
\frac{\delta_0}{2}+\Delta\nonumber
\end{split}
\end{equation}
Function $\lambda(\hat{\theta}_0(\bfx,t_1,\delta))$ does not
belong to $\mathcal{U}_\epsilon([\theta]_\Delta)$, therefore there
exists $\tau(t)$ such that
$|f(\bfx(t+\tau(t),\theta)-f(\bfx(t+\tau(t)),\hat{\theta})|>2\delta_0+\Delta$,
and, consequently
$|\alpha(t+\tau(t))+\varphi(\psi(\bfx(t+\tau(t))))|>\delta_0+\delta_0/2+\Delta$,
$t+\tau(t)>t_1$. Let us denote $t_6=t+\tau(t)$. Notice that
$\dot{\alpha}(t), \pd \varphi(\psi((\bfx(t)))/\pd \psi \dot{\psi}
\in L_\infty$. Hence, there exist constant $\tau_1>0$ such that
$|\varphi(\psi(t))+\alpha(t)|>\delta=\delta_0+\Delta$ for any
$t\in[t_6,t_6+\tau_1]$. Furthermore, function
$\hat{\theta}_0(\bfx,t,\delta)$ is monotonic and, according to
(\ref{main:d_theta_0_mon}), (\ref{nonlinear_PE}) increases by at
least $N \tau_1 \delta^2/2 $ over time interval $[t_1, t]$, $t>
t_1 + N T_1$ until $\lambda(\hat{\theta}_0(t))\in
\mathcal{U}_\epsilon([\theta]_\Delta)$. On the over hand, we have
shown that function $\hat{\theta}_0(\bfx,t,\delta)$ is bounded
from above and, furthermore, converges to the limit
$\hat{\theta}_{0,\infty}<\infty$. Therefore, there exists time
instant $t_7$,  $t_1<t_7<\infty$ such that
$\lambda(\hat{\theta}_0(\bfx,t,\delta))\in
\mathcal{U}_\epsilon([\theta]_\Delta)$ for any $t>t_7$.

Let us assume that $\lambda(\hat{\theta}_0(\bfx,t_1,\delta))$
$\in\mathcal{U}_\epsilon([\theta]_\Delta)$. Then there are two
possibilities. First, $\lambda(\hat{\theta}_{0}(\bfx,t,\delta))$
$\in \mathcal{U}_\epsilon([\theta]_\Delta)$ for all $t>t_1$, and
second, $\lambda(\hat{\theta}_{0}(\bfx,t,\delta))$ eventually
leaves the domain $\mathcal{U}_\epsilon([\theta]_\Delta)$. In this
case, however, the previously provided arguments apply. {\it The
theorem is proven.}

{\it Proof of Theorem \ref{theorem:nonlinear_multi_dim}.} Let us
denote
$\upsilon(t)=f(\bfx(t),\thetavec)-f(\bfx(t),\etavec(\hat{\theta}))$.
According to Assumption \ref{assume:multi-dim} function
$\upsilon(t)\in L_\infty$ and moreover
$\|\upsilon(t)\|_\infty=\Delta$.  Then $\bfx$ is bounded. Hence,
applying the same arguments as in the proof of Theorem
\ref{theorem:nonlinear_single_dim} we conclude that
$\hat{\theta}$, $\dot{\hat{\theta}}$ are bounded. Therefore,
$\dot{\upsilon}(t)\in L_\infty$. Then Theorem
\ref{theorem:nonlinear_single_dim} applies, which proves the
theorem.

\end{document}